\documentclass[11pt]{article}
\title
{Embedding spanning trees in random graphs}
\author{Michael Krivelevich
\thanks{School of Mathematical Sciences, Raymond and Beverly
Sackler Faculty of Exact Sciences, Tel Aviv University, Tel Aviv,
69978, Israel. Email: krivelev@post.tau.ac.il. Research supported in
part by a USA-Israel BSF grant and by a grant from the Israel
Science Foundation.}
 }

\usepackage{amsmath, amsfonts}

\begin{document}
\bibliographystyle{plain}
\maketitle
\newtheorem{theorem}{Theorem}
\newtheorem{thm}{Theorem}[section]
\newtheorem{claim}{Claim}
\newtheorem{lemma}[thm]{Lemma}
\newtheorem{propos}{Proposition}
\newtheorem{conjecture}{Conjecture}
\newtheorem{problem}{Problem}
\newtheorem{defin}{Definition}
\newtheorem{corol}{Corollary}
\newcommand{\Proof}{\noindent{\bf Proof.}\ \ }
\newcommand{\Remarks}{\noindent{\bf Remarks:}\ \ }
\newcommand{\whp}{{\bf whp}}
\newcommand{\prob}{probability}
\newcommand{\rn}{random}
\newcommand{\rv}{random variable}
\newcommand{\hpg}{hypergraph}
\newcommand{\hpgs}{hypergraphs}
\newcommand{\subhpg}{subhypergraph}
\newcommand{\subhpgs}{subhypergraphs}
\newcommand{\bH}{{\bf H}}
\newcommand{\cH}{{\cal H}}
\newcommand{\cT}{{\cal T}}
\newcommand{\cF}{{\cal F}}
\newcommand{\cD}{{\cal D}}
\newcommand{\cC}{{\cal C}}

\begin{abstract}
We prove that if $T$ is a tree on $n$ vertices with maximum degree
$\Delta$ and the edge probability $p(n)$ satisfies: $np\ge
C\max\{\Delta\log n,n^{\epsilon}\}$ for some positive $\epsilon>0$,
then with high probability the random graph $G(n,p)$ contains a copy
of $T$. The obtained bound on the edge probability is shown to be
essentially tight for $\Delta=n^{\Theta(1)}$.
\end{abstract}

\section{Introduction}
In this paper we consider the problem of embedding a copy of a given
spanning tree $T$ on $n$ vertices into the binomial random graph
$G(n,p)$.

Embedding problems are one of the most classical subjects in
Extremal and Probabilistic Combinatorics. There is a large variety
of results about finding given subgraphs, or graphs belonging to a
given family, in random graphs. Here we concentrate on embedding
large trees in binomial random graphs.

The problem of embedding large or {\em nearly} spanning trees in
random graphs on $n$ vertices (where by a nearly spanning tree we
mean a tree $T$ whose number of vertices is at most $(1-c)n$ for
some constant $c>0$) is a rather well researched subject, especially
in the case of trees with bounded maximum degree, see, e.g.,
\cite{FdV1}, \cite{AKS}, \cite{FP}, \cite{FdV2}, \cite{H}. In
particular, Alon, Sudakov and the author proved in \cite{AlKS} that
for given $\epsilon>0$ and integer $d$ there exists
$C=C(d,\epsilon)>0$ such that \whp\footnote{An event ${\cal E}_n$
occurs with high probability, or \whp\ for brevity, in the
probability space $G(n,p)$ if
 $\lim_{n\rightarrow\infty} Pr[G\sim G(n,p)\in{\cal E}_n] = 1$.}
 the random graph $G(n,p)$ with $p=C/n$ has a copy of a tree $T$ on
$(1-\epsilon)n$ vertices of maximum degree at most $d$ (in fact
\cite{AlKS} proved that such a random graph contains \whp\ a copy of
{\em every} such tree); better constant dependence and the
resilience version of this result have recently been obtained in
\cite{BCPS} and \cite{BCS}, respectively.

In contrast, nearly nothing has been known for the case of embedding
{\em spanning} trees. Even the case of embedding spanning trees of
bounded maximum degree appears to be unaddressed, apart from some
sporadic cases. Of course, the most classical result is about
embedding a Hamilton path, or even a Hamilton cycle, in $G(n,p)$;
Koml\'os and Szemer\'edi \cite{KS} and independently Bollob\'as
\cite{B} proved that if $p(n)\ge\frac{\ln n+\ln\ln n+\omega(1)}{n}$,
where $\omega(1)$ is any function tending to infinity arbitrarily
slowly with $n$, then \whp\ $G(n,p)$ contains a Hamilton cycle. Alon
et al. \cite{AlKS} observed that if a tree $T$ has a linear in $n$
number of leaves, then \whp\ $G(n,C\ln n/n)$ contains a copy of $T$
for some large enough $C>0$; the proof is not that hard  and
utilizes the embedding result for nearly spanning trees from the
same paper. However, no general result for this problem has been
obtained, and even the case of the {\em comb} (which is the path
$P_0$ of length $\sqrt{n}-1$ with disjoint paths of length
$\sqrt{n}-1$ attached to each vertex of $P_0$), interpolating in
some sense  between the above mentioned solved cases, is open; this
natural question has been communicated to us by Jeff Kahn
\cite{Kahn}.

Here we make a substantial step forward in solving this class of
problems. Our main result if the following embedding theorem.

\begin{theorem}\label{th1}
Let $T$ be a tree on $n$ vertices of maximum degree $\Delta$. Let
$0<\epsilon<1$ be a constant. If
$$
np\ge \frac{40}{\epsilon}\Delta\ln n+n^{\epsilon}\,,
$$
then \whp\ a random graph $G(n,p)$ contains a copy of $T$.
\end{theorem}

In other words, starting from $\Delta(T)=n^{\epsilon}$, edge
probability $p=\frac{C\Delta\ln n}{n}$ is enough to get \whp\ a copy
of $T$ in $G(n,p)$.

It is not hard to see that the dependence of $p$ on $\Delta$, posted
in Theorem \ref{th1}, is optimal up to a constant factor in the
range $np=n^{\Theta(1)}$. In order to state this result formally,
for integers $n\ge\Delta\ge 3$ define the tree $T(n,\Delta)$ as
follows. Write $n=(\Delta-1)k-r$, where $0\le r\le \Delta -2$. Take
a path $P=(v_1,\ldots,v_k)$ with $k$ vertices, attach to
$v_1,\ldots,v_{k-1}$ vertex disjoint stars with $\Delta-2$ leaves
each, and finally attach to $v_k$ a star with $\Delta-2-r$ leaves.
The tree $T(n,\Delta)$ has $n$ vertices and is of maximum degree at
most $\Delta$. For future reference observe that the $k$ vertices of
$P$ dominate the remaining $n-k$ vertices of $T(n,\Delta)$.

\begin{theorem}\label{th2}
For every $\epsilon>0$ there exists $\delta>0$ such that if
$n^{\epsilon}\le \Delta\le \frac{n}{\ln n}$, then a random graph
$G(n,p)$ with $p=\frac{\delta\Delta\ln n}{n}$ \whp\ does not contain
a copy of $T(n,\Delta)$.
\end{theorem}

There is a certain similarity in appearances between the above
results and the theorem of Koml\'os, S\'ark\"ozy and Szemer\'edi
\cite{KSS}, who proved that for $\delta>0$ and all large enough $n$,
{\em any} graph $G$ on $n$ vertices of minimum degree at least
$(1/2+\delta)n$ contains a copy of {\em every} tree $T$ on $n$
vertices of maximum degree $\Delta(T)\le cn/\ln n$, where
$c=c(\delta)$ is a small enough constant; they noticed that their
condition on $\Delta(T)$ is essentially tight too (actually because
of the random graph $G(n,p)$ with $1/2<p<1$ and the above described
tree $T(n,\Delta)$, just like in our Theorem \ref{th2}). The
arguments of \cite{KSS} are naturally very different and do not seem
to have much bearing on the situation in (sparse) random graphs.

In order to ease the reader's task we now give a brief description
of the proof of Theorem \ref{th1}. The key definition used is that
of a bare path:

\begin{defin}
A path $P$ in a tree $T$ is called {\em bare} if all vertices of $P$
have degree exactly two in $T$.
\end{defin}

In the proof of Theorem \ref{th1}, we first argue that every tree
$T$ on $n$ vertices has a linear in $n$ number of leaves, or a
collection of vertex disjoint bare paths of (large) constant length
each (Lemma \ref{le1}). The former case is rather easy; similarly to
the argument outlined in \cite{AlKS}, we first embed the subtree $F$
of $T$, obtained by deleting from $T$ a linear number of leaves, by
a straightforward greedy algorithm (Lemma \ref{le2}). Then we embed
the remaining edges between the omitted leaves of $T$ and their
fathers (Lemma \ref{le3}); the restriction $np\ge C\Delta\ln n$ is
induced by this part. In the complementary case, where the number of
leaves of $T$ is relatively small, we first take out a linear number
of disjoint constant length bare paths to obtain a subforest $F$ of
$T$; we embed $F$ using the same greedy argument (Lemma \ref{le2}).
Then we are left with embedding the remaining bare paths; we do this
by reducing the problem to that of finding a factor of cycles (with
some extra conditions imposed) in a random graph, and then by
invoking a beautiful result of Johansson, Kahn and Vu \cite{JKV}
about factors in random graphs (Lemma \ref{le4}). In both above
described cases we need our random edges to come in two independent
chunks: the standard trick of representing $G\sim G(n,p)$ as
$G=G_1\cup G_2$, where $G_i\sim G(n,p_i)$ and $1-p=(1-p_1)(1-p_2)$,
allows for this readily.

The notation used in the paper is pretty standard. We systematically
suppress rounding signs for the sake of clarity of presentation.

The proofs of Theorems \ref{th1} and \ref{th2} are given in the next
section. The last section of the paper is devoted to concluding
remarks.

\section{Proofs}
\begin{lemma}\label{le1}
Let $k,l,n>0$ be integers. Let $T$ be a tree on $n$ vertices with at
most $l$ leaves. Then $T$ contains a collection of at least
$\frac{n-(2l-2)(k+1)}{k+1}$ vertex disjoint bare paths of length $k$
each.
\end{lemma}

\Proof Define
\begin{eqnarray*}
V_1 &=& \{v\in V(T): d(v)=1\}\,,\\
V_2 &=& \{v\in V(T): d(v)=2\}\,,\\
V_3 &=& \{v\in V(T): d(v)\ge 3\}\,.
\end{eqnarray*}
Clearly $V_1$ is the set of leaves of $T$ and thus satisfies
$|V_1|\le l$. We have:
\begin{eqnarray*}
2n-2 = 2|E(T)| =\sum_{v\in V(T)}d(v) &\ge& |V_1|+2|V_2|+3|V_3|\\
&=& 2(|V_1|+|V_2|+|V_3|)+ (|V_3|-|V_1|)\\
&=& 2n + |V_3|-|V_1|\,,
\end{eqnarray*}
implying $|V_3|\le |V_1|-2\le l-2$.

$T$ has $|V_1|+|V_3|-1\le 2l-3$ internally disjoint paths connecting
between the vertices of $V_1\cup V_3$, with all internal vertices of
these paths being of degree two. In each such path, pick a largest
collection of vertex disjoint subpaths of length $k$. This leaves at
most $k$ vertices of the path uncovered, so altogether the so formed
collection of bare paths of length $k$ in $T$ contains all but at
most $(|V_1|+|V_3|)+(|V_1|+|V_3|-1)k<(2l-2)(k+1)$ vertices, implying
that the total number of paths of length $k$ in the collection is at
least $\frac{n-(2l-2)(k+1)}{k+1}$ as required. \hfill $\Box$

\begin{lemma}\label{le2}
Let $0<a<1$ be a constant. Let $F$ be a tree on $(1-a)n$ vertices of
maximum degree $\Delta$. If $anp\ge 3\Delta+5\ln n$, then \whp\ a
random graph $G(n,p)$ contains a copy of $F$.
\end{lemma}

\Proof Choose arbitrarily a root $r$ of $F$ and fix some search
order $\pi$, say BFS, on $F$ starting from $r$. Let
$\pi=(v_1=r,\ldots,v_m)$ with $m=(1-a)n$. We will embed $F$ in
$G\sim G(n,p)$ according to $\pi$. Let $\phi$ be the so constructed
embedding.

Suppose we are to embed the children of a current vertex $v_i,\ 1\le
i\le m-1$, in $G$. Let $U_i\subset [n]$ be the set of vertices
already used for embedding, clearly $|U_i|<m$. Expose the edges of
$G$ from $\phi(v_i)$ to $[n]-U_i$. We need to find at most
$d_F(v_i)\le \Delta$ neighbors of $\phi(v_i)$ outside $U_i$. The
probability of this not happening is at most
$$
Pr[Bin(n-m,p)<\Delta] \le
e^{-\frac{(anp-\Delta)^2}{2anp}}<e^{-\frac{2anp}{9}} \ll
\frac{1}{n}\ .
$$
Taking the union bound over all embedding steps, we conclude that
\whp\ $G$ contains a copy of $F$. \hfill$\Box$

\begin{lemma}\label{le3}
Let $0<d_1,\ldots,d_k$ be integers satisfying: $d_i\le \Delta$,
$\sum_{i=1}^kd_i=l$. Let $A=\{a_1,\ldots,a_k\}$, $B$ be disjoint
sets of vertices with $|B|=l$. Let $G$ be a random bipartite graph
with sides $A$ and $B$, where each pair $(a,b)$, $a\in A,b\in B$, is
an edge of $G$ with probability $p$, independently of other pairs.
If
$$
p\ge \frac{2\Delta\ln l}{l}\,,
$$
then \whp\ as $l\rightarrow\infty$ the random graph $G$ contains a
collection $S_1,\ldots, S_k$ of vertex disjoint stars such that
$S_i$ is centered at $a_i$ and has the remaining $d_i$ vertices in
$B$.
\end{lemma}

\Proof Define an auxiliary (random) bipartite graph $G'$ with sides
$A'$ and $B$, where $|A'|=|B|=l$. The vertices of $A'$ are
partitioned into $k$ pairwise disjoint sets $A_1,\ldots, A_k$ with
$|A_i|=d_i,\ 1\le i\le k$. $G$ has an edge between $a'\in A_i$ and
$b\in B$ with probability $p_i$, where $(1-p_i)^{d_i}=1-p$, implying
$p_i\ge p/d_i\ge p/\Delta$. The distribution $G'$ induces the
distribution of $G$ by the obvious projection: $G$ has an edge
between $a\in A$ and $b\in B$ iff $G'$ has some edge between $A_i$
and $B$. Observe that if $G'$ has a perfect matching $M'$ then $G$
has the desired collection of stars $\{S_i\}$, obtained by
projecting $A'$ back into $A$ (the vertices of $A_i$ are projected
onto $a_i$).

By the classical results about random graphs (see, e.g., Section 7.3
of \cite{Bol}) and the monotonicity of the property of having a
perfect matching it is enough to require that all individual edge
probabilities in $G'$ are at least $(\ln l+\omega(1))/l$. Recalling
that $p_i\ge p/\Delta$, we see that the lemma's assumption $p\ge
\frac{2\Delta\ln l}{l}$ guarantees the required condition. \hfill
$\Box$

\begin{lemma}\label{le4}
Let $k\ge 3$ be a fixed integer. Let $G$ be distributed as
$G((k+1)n_0,p)$. Let $S=\{s_1,\ldots,s_{n_0}\}$,
$T=\{t_1,\ldots,t_{n_0}\}$ be disjoint vertex subsets of
$[(k+1)n_0]$. If
$$
p \ge C\left(\frac{\ln n_0}{n_0^{k-1}}\right)^{1/k}\,,
$$
for some large enough constant $C=C(k)$, then \whp\ $G$ contains a
family $\{P_i\}_{i=1}^{n_0}$ of vertex disjoint paths, where $P_i$
is a path of length $k$ connecting $s_i$ and $t_i$.
\end{lemma}

\Proof Fix a partition of $V(G)-S\cup T$ into vertex disjoint
subsets $V_1,\ldots,V_{k-1}$ of cardinality $|V_i|=n_0$ each. Define
an auxiliary graph $H$ with vertex set $V(H)=X\cup V_1\cup\ldots\cup
V_{k-1}$, where $X=\{x_1,\ldots,x_{n_0}\}$. For $1\le i\le k-2$, the
edges of $H$ between $V_i$ and $V_{i+1}$ are identical to those of
$G$. For $v\in V_1$ and $x_j\in X$, $(v,x_j)$ is an edge of $H$ iff
$(v,s_j)$ is an edge of $G$. Similarly, for $v\in V_{k-1}$ and
$x_j\in X$, $(v,x_j)$ is an edge of $H$ iff $(v,t_j)$ is an edge of
$G$. Notice that each relevant pair in $V(H)$ becomes an edge of $H$
independently and with probability $p$.

Suppose now that $H$ contains a $C_k$-factor
$\{S_1,\ldots,S_{n_0}\}$, where each cycle $S_j$ traverses the sets
$X,V_1,\ldots,V_{k-1}$ in this order. Each such cycle $S_i$
translates to a path of length $k$ between $s_i$ and $t_i$ in $G$,
and these paths are pairwise disjoint.

It thus remains to argue that the random graph $H$ contains \whp\
the desired collection of cycles. This can be obtained from the
result of Johansson, Kahn and Vu \cite{JKV} through straightforward
(but quite tedious) modification of their arguments. (They proved
that a random graph $G(kn,p)$ with $p\ge C(k) \left(\frac{\ln
n}{n^{k-1}}\right)^{1/k}$ contains \whp\ a factor of cycles $C_k$,
we need the factor in a $k$-{\em partite} random graph; moreover,
the cycles in the factor are required to traverse the parts in the
{\em prescribed} order.) \hfill$\Box$

\bigskip

\noindent{\bf Proof of Theorem \ref{th1}.} Set
\begin{eqnarray*}
\delta &=& \frac{\epsilon}{10}\,,\\
k &=& \left\lceil \frac{2}{\epsilon}\right\rceil\,.
\end{eqnarray*}
We consider two cases.

\noindent{\bf Case 1.} $T$ has at least $\delta n$ leaves.\\
We represent $G$ as the union $E(G)=E(G_1)\cup E(G_2)$, where $G_1$,
$G_2$ are two independent random graphs, both distributed according
to $G(n,p')$ with $1-p=(1-p')^2$ (and thus $p'\ge p/2$). Let $F$ be
a subtree of $T$ obtained by deleting from $T$ an arbitrary set of
$\delta n$ leaves. We first find a copy $\phi(F)$ of $F$ in $G_1$ --
such a copy exists \whp\ due to Lemma \ref{le2}. Let now
$B=[n]-V(\phi(F))$, and let $A\subset V(\phi(F))$  be the set of
images of the fathers of the $\delta n$ leaves deleted from $T$ to
form $F$. Denote $A=\{a_1,\ldots,a_k\}$, and let $d_i\le \Delta$ be
the number of leaves in $T$ connected to the preimage
$\phi^{-1}(a_i)$ and left outside $F$; clearly
$\sum_{i=1}^{k}d_i=\delta n$. In order to complete the embedding of
$T$ into $G$, we need to find in $G$ $k$ vertex disjoint stars
$S_1,\ldots,S_k$, where the star $S_i$ is centered in $a_i$ and has
the remaining $d_i$ vertices in $B$. We invoke Lemma \ref{le3} to
find such stars \whp\ using the (random) edges of $G_2$ between $A$
and $B$. Since the edge probability in $G_2$ is at least $p/2$, we
need to verify that
$$
\frac{p}{2}\ge \frac{2\Delta\ln(\delta n)}{\delta n}\ .
$$
Recalling that $\delta=\frac{\epsilon}{10}$ and $p\ge
\frac{40\Delta\ln n}{\epsilon n}$, we see that this condition is
fulfilled indeed.

\noindent{\bf Case 2.} $T$ has less than $\delta n$ leaves.\\
We again represent $G$ as the union $E(G)=E(G_1)\cup E(G_2)$ as in
the previous case. According to Lemma \ref{le1}, $T$ contains a
family of $n_0=\frac{n-(2\delta n-2)(k+1)}{k+1}=\Theta(n)$ of vertex
disjoint bare paths of length $k$. Let $F$ be a subforest of $T$
obtained by deleting the internal vertices of  such a family of bare
paths. We first use the edges of $G_1$ to find \whp\ a copy
$\phi(F)$ of $F$; this is possible again due to Lemma \ref{le2}. It
now remains to insert these $n_0$ bare paths, connecting between
prescribed pairs of vertices. We can apply Lemma \ref{le4} to the
edges of $G_2$ to meet this goal. Since the edge probability in
$G_2$ is at least $p/2$ and
$$
\frac{p}{2}\ge \frac{n^{-1+\epsilon}}{2}\gg \left(\frac{\ln
n}{n^{k-1}}\right)^{1/k}
$$
(recall $\epsilon>1/k$), the graph $G_2$ contains indeed the
required collection of paths \whp. The proof is
complete.\hfill$\Box$

\bigskip

\noindent{\bf Proof of Theorem \ref{th2}.} Set
$$
k = \left\lceil\frac{n}{\Delta-1}\right\rceil\,.
$$
Consider the random graph $G(n,p)$ with $p=\frac{\delta\Delta\ln
n}{n}$, the value of $\delta=\delta(\epsilon)$ to be chosen later.
Recall that in the tree $T(n,\Delta)$ $k$ vertices of the spine path
$P$ dominate the rest of the graph. It thus suffices to show that
\whp\ $G(n,p)$ has no dominating set of size $k$. The probability
that such a dominating set exists is at most
\begin{gather*}
\binom{n}{k}(1-(1-p)^k)^{n-k}\le \left(\frac{en}{k}\right)^k\,
e^{-(n-k)(1-p)^k}\\
\le (3\Delta)^{\frac{n}{\Delta}}\, e^{-\frac{n}{2}e^{-pk}} \le
e^{\frac{2n\ln\Delta}{\Delta}-\frac{n^{1-\delta}}{3}}\\
\le e^{2n^{1-\epsilon}\ln n-\frac{n^{1-\delta}}{3}}\ .
\end{gather*}
Taking $\delta=\epsilon/2$ we see that \whp\ the random graph
$G(n,p)$ does not contain a dominating set of size $k$ and thus
\whp\ does not contain a copy of $T(n,\Delta)$.\hfill$\Box$

\section{Concluding remarks}
We have shown that the (pretty immediate) lower bound on the edge
probability $p(n)\ge c\Delta(T)\ln n/n$ for the random graph
$G(n,p)$ to contain \whp\ a copy of a given spanning tree $T$ of
maximum degree $\Delta$ is tight up to a constant factor in the
range $\Delta=n^{\Theta(1)}$.

The regime $\Delta(T)=n^{o(1)}$ stays largely open. In particular,
we were not able to provide a satisfactory solution for the most
natural case of embedding spanning trees with bounded maximum
degree. Our result only shows that in this case is is enough to
require $p(n)=n^{-1+o(1)}$; this is probably not the tightest bound
possible.

For the case of embedding a bounded degree spanning tree $T$ with
$cn$ leaves \cite{AlKS} has shown that it is enough to take
$p(n)=C\ln n/n$, where $C$ may depend on $c$. It is unclear whether
such a dependence is necessary. It seems plausible that assuming
$p(n)=(1+o(1))\ln n/n$ may be enough.

Finally, it would be very interesting to obtain sufficient
conditions for embedding spanning trees with given maximum degree
applicable to {\em pseudo-random} graphs.

\end{document}